\documentclass[12pt]{amsart}
\usepackage{amssymb,amsthm,amsmath,a4}
\newtheorem{thm}{Theorem}[]
\newtheorem{theorem}[thm]{Theorem}

\newtheorem{lemma}[thm]{Lemma}

\theoremstyle{definition}

\newtheorem{example}[thm]{Example}

\theoremstyle{remark}
\newtheorem{remark}[thm]{Remark}
\newcommand{\C}{\mathbb C}
\newcommand{\N}{\mathbb N}
\newcommand{\dr}{\mathrm{d}}
\newcommand{\DC}{\mathcal{D}}
\newcommand{\EC}{\mathcal{E}}
\newcommand{\GC}{\mathcal{G}}
\newcommand{\HC}{\mathcal{H}}
\newcommand{\QC}{\mathcal{Q}}
\newcommand{\If}{{\mathfrak I}}
\newcommand{\LC}{\mathcal{L}}

\newcommand{\R}{\mathbb R}
\newcommand{\Sf}{\mathfrak S}
\newcommand{\vb}{\mathbf{v}}
\newcommand{\xb}{\mathbf{x}}

\newcommand{\de}{\delta}
\newcommand{\eps}{\varepsilon}
\newcommand{\varf}{\varphi}
\newcommand{\la}{\lambda}
\newcommand{\La}{\Lambda}

\newcommand{\Om}{\Omega}
\newcommand{\si}{\sigma}

\renewcommand{\th}{\theta}
\newcommand{\X}{\Xi}

\newcommand{\rank}{\mathrm{rank}\,}
\newcommand{\re}{\mathrm{Re\,}}
\newcommand{\im}{\mathrm{Im\,}}
\newcommand{\Tr}{\mathrm{Tr}\,}

\def\leq{\leqslant}
\def\geq{\geqslant}

\begin{document}

\title[The Berezin and G\aa rding inequalities]
{The Berezin and G\aa rding inequalities}
\author{Y. Safarov}
\address{Department of Mathematics, King's College,
Strand, London WC2R 2LS}
\email{ysafarov@mth.kcl.ac.uk}
\date{December 2004}


\begin{abstract}
Let $\,\varf\,$ be a convex function on $\,\mathbb C\,$,
$\,\LC(\si)\,$ be a pseudodifferential operator with symbol
$\,\si\,$, $\,\La_\si\,$ be the set of its eigenvalues and
$\,m(\la)\,$ be the multiplicity of an eigenvalue
$\,\la\in\La_\si\,$. Under certain natural assumptions about
properties of pseudodifferential operators, we prove that
$\,\sum_{\la\in\La_\si}m(\la)\,\varf(\la)\le\re
\Tr\LC(\varf(\si))+R\,$, where $\,R\,$ is an error term of the same
order as the remainder term in the G\aa rding inequality.
\end{abstract}

\maketitle

\section{Introduction}\label{S1}

Let $\varf:\R\to\R$ be a convex function, $B$ be a self-adjoint
operator and $P$ be an orthogonal projection in a separable
Hilbert space $H$. Then
\begin{equation}\label{1}
\Tr\varf(P\left.B\right|_{PH})\ \leq\
\Tr(P\left.\varf(B)\right|_{PH})\,,
\end{equation}
provided that the operators in the right and left hand sides are
well defined and belong to the trace class $\Sf_1$. This estimate
was proved in \cite{B1} and is often called ``the Berezin
inequality'' (some generalizations of (1) were obtained in
\cite{LS1}). By the spectral theorem, every self-adjoint operator is
unitary equivalent to a multiplication operator in $L_2\,$.
Therefore Berezin's result can be reformulated in the following way:
if $\varf:\R\to\R$ is a convex function, $b$ is a real-valued
function, $\{b\}$ is the corresponding multiplication operator,
$U:H\to L_2$ is an isometry onto a subspace of $L_2$ and
$\QC(b):=U^*\{b\}U$, then
\begin{equation}\label{2}
\Tr\varf(\QC(b)))\ \leq\ \Tr\QC(\varf(b))
\end{equation}
whenever  the operators $\QC(\si)$ and $\QC(\varf(b))$ are well
defined and belong to $\Sf_1$.

The Berezin inequality has been used for the study of spectral
properties of differential and pseudodifferential operators. If $B$
is a self-adjoint pseudodifferential operator with symbol $\si_B$
then, under certain assumptions, $\varf(B)$ is a pseudodifferential
operator whose symbol coincides with $\varf(\si_B)$ modulo a lower
order term. In this case the right hand side of (\ref{1}) is equal
to the sum of an integral of $\varf(\si_B)$ and a lower order
remainder, and (1) implies asymptotic formulas for the spectrum of
the operator $P\left.B\right|_{PH}$ (see, for example, \cite{LS2}).
If there exist an isometry $U:H\to L_2$ and a function $b$ such that
$B=\QC(b$), and $\Tr\QC(\varf(b))$ is given by an explicit formula
then (2) yields estimates for the spectrum of the operator $B$
itself (see, for instance, \cite{ELSS}).

The first scheme works only for self-adjoint pseudodifferential
operators $B$ and relies on symbolic functional calculus. The second
allows one to obtain estimates only in terms of the function $b$
which depends on the choice of the isometry $U$. The main problem in
this scheme is to construct a suitable isometry $U$ and to
investigate the relation between $b$ and the actual symbol $\si_B$.
For some operators this can be done with the use of the so-called
coherent states (as in \cite{B2} or \cite{ELSS}). With the exception
of some very special cases, the formulas relating $b$ and $\si_B$
contain a lower order error term. These formulas, together with (2),
imply asymptotic estimates in terms of $\si_B$ with a similar error
term.

The aim of this paper is to show that the Berezin inequality is an
elementary consequence of the G\aa rding inequality. If the G\aa
rding inequality holds with a lower order error term then the
Berezin inequality contains an error term of the same order. Note
that the G\aa rding inequality is a simpler result than a symbolic
functional calculus or a coherent state representation and, as a
rule, immediately follows from either of these two.

The inequalities (\ref{1}) and (\ref{2}) are easily proved by
representing the quadratic form of $P\left.B\right|_{PH}$ or $B$
as a Lebesgue integral and applying Jensen's inequality. Our proof
does not involve the spectral theorem or Lebesgue integrals.
Instead, we observe that Jensen's inequality holds for more
general functionals and apply it to suitably chosen functionals on
the space of symbols.

\section{Convex functions}\label{S2}

In this section we shall briefly recall some results from convex
analysis. Keeping in mind possible applications to operator-valued
functions $\si$, we shall consider convex functions on an infinite
dimensional locally convex real vector space $X$. All results and
their proofs are elementary and cannot be substantially simplified
even if $X=\R$.

A function $\varf:X\to[-\infty,+\infty]$ is called {\sl convex} if
its epigraph $\EC(\varf):=\{(t,\xb)\in \R\times
X:t\geq\varf(\xb)\}$ is a convex subset of $\R\times X$. A convex
function is said to be {\sl proper\/} if $|\varf|\not\equiv\infty$
and {\sl closed\/} (or lower semicontinuous) if $\,\EC(\varf)\,$
is closed in $\R\times X$ in the product topology. If $\EC(\varf)$
lies to one side of a hyperplane $\HC$ passing through the point
$(\varf(\xb),\xb)$ then $\HC$ is said to be a {\sl supporting}
hyperplane at $(\varf(\xb),\xb)$. If $\varf$ is convex and $\dim
X<\infty$ then each point $(\varf(\xb),\xb)$ has at least one
(possibly, vertical) supporting hyperplane. In the infinite
dimensional case there may be no supporting hyperplanes.

\begin{example}\label{E1}
Let $X=\R^\infty$ be the space of real sequences
$\xb=\{x_1,x_2,\ldots\}$ provided with the topology of
element-wise convergence, $\,\varf(\xb):=\sum_{j=1}^\infty
x_j^{-1}$ if $x_j>0$ for all $j$ and the sum is finite, and
$\varf(\xb):=+\infty$ otherwise. Then the closed convex set
$\EC(\varf)$ does not have any supporting hyperplanes; in other
words, no linear continuous functional attains its minimal value
on $\EC(\varf)$.
\end{example}

If the supporting hyperplane $\HC$ is not vertical then it
coincides with the graph of an affine function
$\,l_{t^*,\xb^*}(\xb):=\langle\xb^*,\xb\rangle-t^*\,$ with some
$t^*\in\R$ and $\xb^*$ from the dual space $X^*$. The set of
vectors $\xb^*\in X^*$ generating non-vertical supporting
hyperplanes at $(\varf(\xb),\xb)$ is called the {\sl
subdifferential} of $\varf$ at the point $\xb$ and is denoted by
$\partial\varf(\xb)$.

The closed convex function $\,\varf^*(\xb^*):= \sup_{\xb\in
X}\,\{\langle \xb^*,\xb\rangle-\varf(\xb)\}\,$ on $X^*$ is called
the {\sl conjugate} of $\varf$. The following well known result
(see, for example, \cite{ET} or \cite{G}) is a simple consequence of
the separation theorem.

\begin{lemma}\label{L2}
If $\varf$ is a proper closed convex function then $\varf\equiv
\varf^{**}$.
\end{lemma}

We have $\,(t^*,\xb^*)\in\EC(\varf^*)\,$ if and only if
$\,l_{t^*,\xb^*}(\xb)\leq\varf(\xb)\,$ for all $\xb\in X\,$. Given
$\eps\geq0$ and $\xb_0\in X$, let us denote by
$\partial_\eps\varf(\xb_0)$ the set of points
$(t^*,\xb^*)\in\EC(\varf^*)$ such that $\,\varf(\xb_0)-\eps\leq
l_{t^*,\xb^*}(\xb_0)\,$. In particular,
$$
\partial_0\varf(\xb_0)=\{(t^*,\xb^*)\in\EC(\varf^*)\,:\,
\xb^*\in\partial\varf(\xb_0)\,,\,t^*
=\langle\xb^*,\xb_0\rangle-\varf(\xb_0)\}\,.
$$
The set $\,\partial_0\varf(\xb_0)\,$ may well be empty even in the
case $\,X=\R\,$; if $\,\dim X=\infty\,$ then it may happen that
$\,\partial_0\varf(\xb_0)=\emptyset\,$ for all $\xb_0\in X\,$ (see
Example~\ref{E1}). However, by Lemma~\ref{L2},
$\,\partial_\eps\varf(\xb_0)\ne\emptyset\,$ for each $\eps>0\,$
provided that $\,\varf(\xb_0)<+\infty\,$.

\section{Jensen's inequality}\label{S3}

Let $L_X$ be a set of functions $\si:\X\to X$ defined on a
nonempty set $\,\X\,$, $L_\R$ be a set of real-valued functions on
$\X\,$, $\,\If_X$ be a map from $L_X$ to $X$ and
$\,\If_\R:L_\R\to[-\infty,+\infty]\,$ be a real functional on
$L_\R$.

\begin{lemma}\label{L3}
Let $\varf$ be a proper closed convex function such that
$\,\varf(\si)\in L_\R$. Assume that for each $\eps>0$ there exists a
point $(t^*,\xb^*)\in\partial_\eps\varf(\If_X(\si))$ such that
\begin{enumerate}
\item[(\bf a$_1$)]
$\,l_{t^*,\xb^*}(\si)\in L_\R$ and
$\,l_{t^*,\xb^*}\left(\If_X(\si)\right)
\leq\If_\R\left(l_{t^*,\xb^*}(\si)\right)+C_1$,
\item[(\bf a$_2$)]
$\If_\R\left(l_{t^*,\xb^*}(\si)\right)\leq\If_\R(\varf(\si))+C_2$,
\end{enumerate}
where $C_1$ and $C_2$ are real constants. Then
$\,\varf(\If_X(\si))\leq\If_\R(\varf(\si))+C_1+C_2\,$.
\end{lemma}

\begin{proof}
Since $(t^*,\xb^*)\in\partial_\eps\varf(\If_X(\si))$, the
conditions ({\bf a$_1$}) and ({\bf a$_2$}) imply that
$\,\varf(\If_X(\si))-\eps\leq l_{t^*,\xb^*}\left(\If_X(\si)\right)
\leq\If_\R\left(l_{t^*,\xb^*}(\si)\right)+C_1\leq\If_\R(\varf(\si))+C_1+C_2\,$.
Letting $\eps\to0$, we obtain the required inequality.
\end{proof}

The inequality ({\bf a$_1$}) holds with $\,C_1=0\,$ provided that
\begin{enumerate}
\item[(\bf a$'_1$)]
the functional $\,\If_\R\,$ is linear, $\,\If_\R(1)=1\,$ and
$\If_\R\left(\langle\xb^*,\si\rangle\right)=\langle\xb^*,\If_X(\si)\rangle$.
\end{enumerate}
The condition ({\bf a$_2$}) is fulfilled with $C_2=0$ for all
monotone functionals $\,\If_\R\,$, that is, the functionals
$\,\If_\R\,$ satisfying
\begin{enumerate}
\item[(\bf a$'_2$)]
$\,\If_\R(\si_1)\leq\If_\R(\si_2)\,$ whenever $\si_1,\si_2\in L_\R$
and $\,\si_1(\th)\leq\si_2(\th)\,$ for all $\,\th\in\X\,$.
\end{enumerate}

If $\If_\R$ is the Lebesgue integral with respect to a probability
measure or the normalized Perron integral and $\If_X$ is the
corresponding vector-valued integral understood in the weak sense
then ({\bf a$'_1$}) and ({\bf a$'_2$}) hold for all integrable
functions and all $\,(t^*,\xb^*)\in\R\times X^*\,$. Therefore
Lemma~\ref{L3} implies the standard Jensen's inequality. The
following example is less obvious.

\begin{example}\label{E4}
Let $A$ be a liner positive operator in $L_2(\X)\,$, $L_\R$ be the
space of measurable bounded functions on $\X\,$ and $\,V_\si\,$ be
the operator of multiplication by the function $\si$. Denote by
$\la_1(\si),\la_2(\si),\ldots$ the ordered eigenvalues of the
Friedrichs extension of the operator $A+V_\si$ lying below its
essential spectrum.

Let us fix $n\in\N$ and take $X=\R$ and $\If_\R=\la_n(\si)$. In
view of the Rayleigh--Ritz formula, $\If_\R$ satisfies ({\bf
a}$'_2$) for all $\si\in L_\R$. We have
$$
A+t\,V_\si\ \geq\ \de\,(A+V_\si)+(t-\de)\,\la_n(\si)
+\inf_{\th\in\X}\left((t-\de)\,(\si(\th)-\la_n)\right)
$$
for all $\,\de\in[0,1]\,$ and $\,t\in\R\,$. Let $\varf$ be a
proper closed convex function such that $\varf(\si)\in L_\R\,$.
The above inequality implies ({\bf a}$_1$) with
$$
C_1\ =\ F(\la_n(\si))\ :=\ \max\left\{a_-(\sup\si-\la_n(\si))_+\,,
\,(b-1)_+(\la_n(\si)-\inf\si)\right\},
$$
where $\,a=\sup\partial\varf(\inf\si)\,$,
$\,b=\inf\partial\varf(\sup\si)\,$, and the subscripts $\pm$ denote
the positive and negative parts (we define
$\,\sup\partial\varf(t):=-\infty\,$ and
$\,\inf\partial\varf(t):=+\infty\,$ when
$\partial\varf(t)=\emptyset\,$). Therefore, by Lemma~\ref{L3},
$$
\varf\left(\la_n(\si)\right)
\leq\la_n\left(\varf(\si)\right)+F(\la_n(\si))\,.
$$
Note that $\,F(\la_n(\si))=0\,$ whenever $\,[a,b]\subset[0,1]\,$
or $\,[a,b]\subset[-\infty,1]\,$ and $\,\sup\si\leq\la_n(\si)\,$.
\end{example}

\section{Berezin inequality}\label{S4}

Let $L_\C$ be a linear space of complex-valued functions on $\X\,$
containing the constant functions and closed with respect to the
complex conjugation, and let $L_\R$ be the subspace of real-valued
functions. Consider a linear map $\QC$ from $L_\C$ into the space
of linear operators in a separable Hilbert space $H$ (a
quantization) such that $\QC(1)=I$ and
$\,\DC(\QC(\si))=\DC(\QC(\re\si))\bigcap\DC(\QC(\im\si))\,$. If
$\,\si\in L_\C\,$, let $\Om_\si$ be the numerical range of the
operator $\QC(\si)$, $\La_\si$ be the set of its eigenvalues and
$m(\la)$ be the algebraic multiplicity of the eigenvalue
$\,\la\in\La_\si\,$.

Let us fix a bounded operator $T$ and, given $\nu\in\R$, denote by
$G_\nu$ the set of functions $\si\in L_\R$ such that
\begin{enumerate}
\item[(\bf G)]
$\,\re(\QC(\si)u,u)_H\geq-\,\nu\,(Tu,u)_H\,$ for all $u\in
\DC(\QC(\si))$.
\end{enumerate}

\begin{lemma}\label{L5}
Let $\si\in L_\C$ and $\varf$ be a proper closed convex function
on $\C\,$. Assume that for each $\eps>0$ and each $z\in\Om_\si$
there exists $(t^*,z^*)\in\partial_\eps\varf(z)$ satisfying the
following two conditions:
\begin{enumerate}
\item[(\bf b$_1$)]
$\,\im(\QC(\im(z^*\bar\si))u,u)_H\leq\nu_1(Tu,u)_H\,$ for all
$u\in\DC(\QC(\si))$,
\item[(\bf b$_2$)]
$\varf(\si)-\re(z^*\bar\si)+t^*\in G_{\nu_2}\,$,
\end{enumerate}
where $\nu_1$ and $\nu_2$ are real constants. Then
\begin{equation}\label{3}
\varf\left((\QC(\si)u,u)_H\right)\leq\re(\QC(\varf(\si))u,u)_H+
(\nu_1+\nu_2)\,(Tu,u)_H
\end{equation}
whenever $\,u\in\DC(\QC(\si))\bigcap\DC(\QC(\varf(\si)))$ and
$\|u\|_H=1$.
\end{lemma}

\begin{proof}
Let us identify $\C$ with $\R^2$ so that $\langle
z^*,z\rangle=\re(z^*\bar z)$, $\forall z,z^*\in\C$. Then (\ref{3})
immediately follows from Lemma~\ref{L3} with
$\If_\C(\si)=(\QC(\si)u,u)_H$, $\If_\R(\si)=\re(\QC(\si)u,u)_H$
and $C_j=\nu_j\,(Tu,u)_H\,$.
\end{proof}

\begin{lemma}\label{L6}
Let $\si\in L_\R$ and $\varf$ be a proper closed convex function
on $\R\,$. If
\begin{enumerate}
\item[(\bf b$'_1$)]
the operator $\,\QC(\si)\,$ is symmetric
\end{enumerate} and
for each $\eps>0$ and $z\in\Om_\si$ there exists
$(t^*,z^*)\in\partial_\eps\varf(z)$ satisfying {\rm ({\bf b$_2$})}
then {\rm (\ref{3})} holds with $\nu_1=0$.
\end{lemma}

\begin{proof}
This is a particular case of Lemma~\ref{L5} with the convex
function on $\C$ which is equal to $\varf$ on $\R$ and to
$+\infty$ on $\C\setminus\R$.
\end{proof}

Note that the second condition in Lemmas~\ref{L5} and \ref{L6} is
satisfied whenever $\varf$ is differentiable and
\begin{enumerate}
\item[(\bf b$'_2$)]
$\ \varf(\si)-\varf(z)-\re(\varf'_z(z)\,(\si-z))\in G_{\nu_2}$ for
all $z\in\Om_\si\,$.
\end{enumerate}

\begin{theorem}\label{T7}
Assume that $\si$ and $\varf$ satisfy the conditions of
Lemma~\ref{L5} or Lemma~\ref{L6}, $\QC(\varf(\si))\in\Sf_1\,$,
$T\in\Sf_1$, and at least one of the following two conditions is
fulfilled:
\begin{enumerate}
\item[(\bf c$_1$)] $\varf$ is nonnegative;
\item[(\bf c$_2$)] the set of the generalized
eigenvectors of the operator $\QC(\si)$ is complete.
\end{enumerate}
Then the set
$\,\La_{\si,\varf}^+:=\{\la\in\La_\si\,:\,\varf(\la)>0\}\,$ is
countable, each eigenvalue $\la\in\La_{\si,\varf}^+\,$ has a
finite algebraic multiplicity,
$\,\sum_{\la\in\La_{\si,\varf}^+}m(\la)\,\varf(\la)<\infty\,$ and
\begin{equation}\label{4}
\sum_{\la\in\La_\si}m(\la)\,\varf(\la)\ \leq\
\re\Tr\QC(\varf(\si))+(\nu_1+\nu_2)\,\Tr T\,.
\end{equation}
\end{theorem}

\begin{proof}
Let $\,\vb=\{v_1,v_2,\ldots,v_k\}\,$ be a finite collection of the
generalized eigenvectors of the operator $\,\QC(\si)\,$
corresponding to the eigenvalues $\la_1,\la_2,\ldots,\la_k\,$.
Denote by $H_\vb$ the finite dimensional invariant subspace of
$\,\QC(\si)\,$ spanned by the vectors $\,(\QC(\si)-\la_j)^nv_j\,$,
$j=1,\ldots,k$, $n=0,1,\ldots$ The restriction
$\,\left.\QC(\si)\right|_{H_\vb}\,$ has the same eigenvalues
$\la_1,\la_2,\ldots,\la_k$ whose algebraic multiplicities
$m_\vb(\la_j)$ do not exceed $m(\la_j)\,$. By Schur's lemma, the
operator $\,\left.\QC(\si)\right|_{H_\vb}\,$ is represented by a
triangular matrix with the diagonal entries $\la_j$ in some
orthonormal basis $\{u_j\}\subset H_\vb$. Applying (\ref{3}) to the
vectors $u_j$, we see that
\begin{multline}\label{5}
\sum_{j=1}^km_\vb(\la_j)\,\varf(\la_j)\ \leq\
\re\Tr\Pi_\vb\QC(\varf(\si))+(\nu_1+\nu_2)\,\Tr\Pi_\vb T\\ \leq\
\|\QC(\varf(\si))\|_{\Sf_1}+\|(\nu_1+\nu_2)\,T\|_{\Sf_1}\,,
\end{multline}
where $\Pi_\vb$ is the orthogonal projection onto the subspace
$H_\vb\,$.

If the set $\,\La_{\si,\varf}^+\,$ were uncountable or there were an
eigenvalue $\la\in\La_{\si,\varf}^+\,$ of infinite algebraic
multiplicity then we could find a positive constant $\de$ and an
arbitrarily large collection $\vb$ of eigenvectors $v_j$ such that
$\la_j\in\La_{\si,\varf}^+\,$ and $\,\varf(\la_j)\geq\de\,$. This
contradicts to (\ref{5}). Therefore $\,\La_{\si,\varf}^+\,$ is
countable and $m(\la)<\infty$ for all $\la\in\La_{\si,\varf}^+\,$.
Choosing a sequence of expanding finite sets
$\vb_1\subset\vb_2\subset\vb_3\subset\ldots$ whose union
$\,\bigcup_{n=1}^\infty\vb_n\,$ contains all eigenvectors
corresponding to the eigenvalues $\,\la\in\La_{\si,\varf}^+\,$,
applying (\ref{5}) with $\,\vb=\vb_n\,$ and letting
$\,n\to\infty\,$, we see that
$\,\sum_{\la\in\La_{\si,\varf}^+}m(\la)\,\varf(\la)\,$ is estimated
by the same sum of the trace norms.

In order to prove (\ref{4}), we note that the left hand side of
(\ref{4}) is $\,-\infty\,$ whenever the set
$\,\La_{\si,\varf}^-:=\{\la\in\La_\si\,:\,\varf(\la)<0\}\,$ is
uncountable or contains an eigenvalue of infinite algebraic
multiplicity. Therefore we can assume without loss of generality
that the set of generalized eigenvectors corresponding to the
eigenvalues $\,\la\in\La_{\si,\varf}^+\bigcup\La_{\si,\varf}^-\,$ is
countable. Let us choose a sequence of finite sets
$\vb_1\subset\vb_2\subset\vb_3\subset\ldots$ such that
$\,\bigcup_{n=1}^\infty\vb_n\,$ contains all these eigenvectors and
$\,\overline{\bigcup_{n=1}^\infty H_{\vb_n}}=H'\,$, where $H'$
denotes the closed linear span of all generalized eigenvectors of
the operator $\QC(\si)\,$.

If ({\bf c}$_2$) holds then $\,H'=H\,$. Therefore, taking
$\vb=\vb_n$ in the first inequality (\ref{5}) and letting
$n\to\infty$, we arrive at (\ref{4}). If $\,H'\ne H\,$ and
$\,\varf\geq0\,$ then we choose an orthonormal basis $\{\tilde
u_i\}$ in $(H')^\bot$, apply (\ref{3}) to $\tilde u_i$ and add up
the obtained inequalities and the first inequality (\ref{5}) with
$\vb=\vb_n\,$. Since $\,\varf((\QC(\si)\tilde u_i,\tilde
u_i)_H)\geq0\,$, now (\ref{4}) is proved by letting
$\,n\to\infty\,$.
\end{proof}

If a convex function $\,\varf:\R\to[-\infty,+\infty]\,$ takes
negative values then the set of its zeros consists of at most two
points. In this case the inclusion $\varf(\QC(\si))\in\Sf_1$
implies ({\bf c}$_2$) for each self-adjoint operator
$\,\QC(\si)\,$. Therefore (\ref{2}) is a particular case of
Theorem~\ref{T7} with $T=0$.

\section{Pseudodifferential operators}\label{S5}

In the theory of pseudodifferential operators, $\X$ is the cotangent
bundle $T^*M$ over a domain $M\subset\R^n$ or a manifold $M$ and
$H=L_2(M)$. For $M\subset\R^n$, quantization is defined by the
formula
\begin{equation}\label{6}
\QC_\tau(\si)u(x)=(2\pi h)^{-n}\int_{\R^n}\int_M
e^{ih^{-1}(x-y)\cdot\xi}\si(\tau x+(1-\tau)y,\xi)\,u(y)\,\dr
y\,\dr\xi\,,
\end{equation}
where $\tau$ is a fixed number from the interval $[0,1]$ and $h$ is
a real parameter (see, for example, \cite{DS}). This definition can
be extended to manifolds $M$ (see \cite{Sa}). The functions $\si$ on
$T^*M$ are called $\tau$-symbols; for $\tau=\frac12$ they called
Weyl symbols. In the classical theory of pseudodifferential
operators one takes $h=1$ and defines the order of a symbol $\si$ in
terms of its behaviour for large $\xi\,$. In the semiclassical
theory the order is defined in terms of asymptotic behaviour as
$h\to0$.

If $M\subset\R^n$ then, obviously,
$\,(\QC_\tau(\si)u,u)_H=(u,\QC_{1-\tau}(\bar\si)u)_H\,$. This
equality remains true on a manifold $M$ if the $\tau$-quantization
is defined as in \cite{Sa}. It implies that the estimate ({\bf
b$_1$}) holds for $\,\QC=\QC_{1/2}\,$ with $\nu_1=0$.

If $\si\geq0$ and $T$ is a lower order operator then ({\bf G}) is
known as the G\aa rding inequality or the Fefferman--Phong
inequality (the latter gives a more precise result in terms of the
order of $T$). The constant $\nu$ in this inequality can usually be
estimated by a functional $\GC(\si)$ which involves partial
derivatives of the symbol $\si$ up to a certain order (see, for
example, \cite{DS}, \cite{H}, \cite{T} or \cite{LM}). The condition
({\bf b$_2$}) means that ({\bf G}) holds uniformly on the set of
nonnegative functions $\varf(\si)-\re(z^*\bar\si)+t^*$, which is the
case whenever the functional $\GC$ is uniformly bounded on this set.
One can obtain explicit formulas for $\GC$ by analyzing the known
proofs of the G\aa rding inequality. However, such analysis lies
outside the scope of this paper. Instead, we conclude by giving
three examples which demonstrate possible applications of
Theorem~\ref{T7} (in the last two examples ({\bf b$_2$}) can be
proved directly).

\begin{example}\label{E8}
Let $M$ be a compact $\,n$-dimensional $\,C^\infty$-manifold, $\,H$
be the space of square integrable half-densities on $M$, $\,S^m$ be
the H\"ormander class of symbols and $\,\Psi^m$ be the corresponding
class of classical pseudodifferential operators in $H$. Consider an
elliptic positive pseudodifferential operator $A\in\Psi^1$ and
denote by $\Pi_\mu$ its spectral projection corresponding to the
interval $[0,\mu)$. It is well known that $\,\rank\Pi_\mu<\infty\,$
and
\begin{equation}\label{7}
\Tr(B\,\Pi_\mu)=
(2\pi)^{-n}\int_{\si_A(x,\xi)<\mu}\si_B(x,\xi)\,\dr
x\,\dr\xi\;+\;O(\mu^{n+m-1})\,,\quad\mu\to\infty\,,
\end{equation}
for every $B\in\Psi^m$ provided that $\,n+m-1\geq0$, where $\si_A$
and $\si_B$ are principal symbols of the operators $A$ and $B$ (see,
for example, \cite{H} or \cite{SV}).

Let $L_\C=S^0$, the quantization $\,\QC_{1/2}:S^0\to\Psi^0\,$ be
defined as in \cite{Sa},
$\,\QC(\si):=\left.\Pi_\mu\QC_{1/2}(\si)\right|_{\Pi_\mu H}\,$ and
$\,T=\left.\Pi_\mu A^{-1}\right|_{\Pi_\mu H}\,$. If $\varf\in
C^\infty(\C)$ is a convex function then, by the G\aa rding
inequality, we have ({\bf b$'_2$}) with some constant $\nu_2$
depending on $\si$ and $\varf$. Since $\,\si-\si_B\in S^{-1}\,$
whenever $\,B=\QC_{1/2}(\si)\,$ (see \cite{Sa}), Theorem~\ref{T7}
and (\ref{7}) imply that
$$
\sum_{\la\in\La_{\mu,B}}m(\la)\,\varf(\la)\ \leq\
(2\pi)^{-n}\int_{\si_A(x,\xi)<\mu}\varf(\si_B(x,\xi))\,\dr
x\,\dr\xi\;+\;O(\mu^{n-1})\,,\quad\mu\to\infty
$$
for every operator $B\in\Psi^0$, where $\,\La_{\mu,B}\,$ is the set
of eigenvalues of $\left.\Pi_\mu B\right|_{\Pi_\mu H}$ and
$\,m(\la)\,$ is the algebraic multiplicity of $\la$.

\end{example}

\begin{remark}\label{R9}
The above inequality was obtained in \cite{LS2} for self-adjoint
operators $\,B\in\Psi^0\,$.
\end{remark}

In the following examples $M\subset\R^n$ is an open bounded set
and $D$ is an open bounded subset of $T^*M\,$.

\begin{example}\label{E10}
Let $\si$ be the characteristic function of $D$, $\,\si':=1-\si\,$,
$L_\C$ be the linear space spanned by $\si$ and $\si'$, and let
$\,\QC=\QC_1\,$. If $\,R:=\QC_1(\si)\,\QC_0(\si')\,$ then
$\,\QC_1(\si)=\QC_1(\si)\,\QC_0(\si)+R\,$,
$\,\QC_1(\si')=\QC_1(\si')\,\QC_0(\si')+R^*\,$ and
$$
2\,\im(\QC_1(\im(z^*\si))u,u)_H=\im
z^*((R-R^*)u,u)_H\,,\qquad\forall z^*\in\C\,.
$$
Therefore the conditions ({\bf b$_1$}) and ({\bf b$'_2$}) are
satisfied with $\,T=|\re R|+|\im R|\,$,
$\,\nu_1=\sup_{z\in\Om_\si}|\im\varf'_z(z)|\,$ and
$\,\nu_2=\sup_{z\in\Om_\si}|\varf(1)-\re\varf'_z(z)|\,$ for every
nonnegative differentiable convex function $\,\varf:\C\to\R\,$
vanishing at the origin $z=0$. The operators $\,\QC_1(\si)\,$ and
$\,\QC_0(\si')\,$ belong to the Hilbert-Schmidt class $\,\Sf_2\,$.
Therefore $\,T\in\Sf_1\,$, $\,\QC_1(\varf(\si))\in\Sf_1\,$, and
Theorem~\ref{T7} implies that
$$
\sum_{\la\in\La_\si}m(\la)\,\varf(\la)\ \leq\ (2\pi
h)^{-n}\,\varf(1)\int_D \dr x\,\dr\xi\;+\;4\left(\varf(1)+
\sup_{z\in\Om_\si}|\varf'_z(z)|\right)\|R\|_{\Sf_1}\,.
$$
Approximating $\si$ and $\si'$ by smooth functions, which are
supported in $D$ and $T^*M\setminus D$ respectively and vanish near
the boundary $\partial D$, we see that $\|R\|_{\Sf_1}=o(h^{-n})$ as
$h\to0$. If the Minkowski dimension $d$ of $\partial D$ is strictly
smaller than $n$ then one can improve this estimate and show that
$\|R\|_{\Sf_1}=O(h^{\alpha-n})$ where $\alpha$ is a positive
constant depending on $d$.
\end{example}

\begin{example}\label{E11}
Let $\,C_0^\kappa(D)\,$ be the subspace of the H\"older space
$\,C^\kappa(\R^{2n})\,$ which consists of real-valued function
vanishing outside $D$. Let us fix $\eps>0$, denote by $L_\R$ the
real linear space spanned by $\,C_0^{n/2+2+\eps}(D)\,$ and constant
functions, and consider the Weyl quantization on $L_\C=L_\R\,$.

Let $\,\si\in C_0^{n/2+2+\eps}(D)\,$ and $\varf$ be a nonnegative
convex function on the closure $\overline{\Om_\si}$ such that
$\,\varf\in C^{n/2+2+\eps}\,$ and $\varf''\geq\de>0$ (recall that
$\,\Om_\si\subset\R)\,$. Then
$\,\varf_z(\si):=\varf(\si)-\varf(z)-\varf'(z)\,(\si-z)=(\psi_z(\si))^2\,$
where $\,\psi_z\in C^{n/2+2+\eps}(\R^{2n})\,$. Expanding the
function $\,\psi_z((x+y)/2,\xi)\,$ by Taylor's formula at $x=y$ and
$y=x$, replacing $(x-y)\,e^{ih^{-1}(x-y)\cdot\xi}$ with
$-ih\,\nabla_\xi e^{ih^{-1}(x-y)\cdot\xi}$ in (\ref{6}) and
integrating by parts with respect to $\xi$, one can easily prove
that
$\,\QC_{1/2}(\varf_z(\si))-\QC_1(\psi_z(\si))\,\QC_0(\psi_z(\si))\,$
coincides with a finite sum of operators $R_k$ whose Schwartz
kernels are given by oscillatory integrals of the form
$$
i\,h^{1-n}\int_{\R^n}\int_0^1\int_0^1 e^{ih^{-1}(x-y)\cdot\xi}
a_{k,1}(l_{k,1}(t_1,x,y),\xi)\,a_{k,2}(l_{k,2}(t_2,x,y),\xi)\,\dr
t_1\,\dr t_2\,\dr\xi\,,
$$
where $\,a_{k,j}\,$ is a derivative of $\psi_z(\si)$ of order $\,0$,
$1$ or $2$ and $\,l_{k,j}(t,x,y)\,$ is one of the following
functions: $x$, $y$, $\,x+t(y-x)/2\,$ or $\,y+t(x-y)/2\,$. The
amplitudes in these oscillatory integrals vanish for all
sufficiently large $\,\xi\,$ because the functions $\,\psi_z(\si)\,$
are constant outside $\,D\,$ and each amplitude contains at least
one derivative of $\,\psi_z\,$ of order 1 or 2. Furthermore, we are
only interested in $\,x,y\in M\,$. Therefore we can replace
$\,a_{k,j}(l_{k,j}(t_j,x,y),\xi)\,$ with
$\,b_{k,j}(l_{k,j}(t_j,x,y),\xi):=
\chi(l_{k,j}(t_j,x,y),\xi)\,a_{k,j}(l_{k,j}(t_j,x,y),\xi)\,$ where
$\,\chi\,$ is a $\,C_0^\infty\,$-function on $\,\R^{2n}\,$ which is
equal to one on a sufficiently large ball. We have
\begin{multline*}
b_{k,1}(l_{k,1}(t_1,x,y),\xi)\,b_{k,2}(l_{k,2}(t_2,x,y),\xi)\\ =\
(2\pi)^{-2n}\int_{\R^n}\int_{\R^n}e^{i\eta_1\cdot
l_{k,1}(t_1,x,y)}\hat b_{k,1}(\eta_1,\xi)\,e^{i\eta_2\cdot
l_{k,2}(t_2,x,y)}\,\hat b_{k,2}(\eta_2,\xi)\,\dr\eta_1\,\dr\eta_2\,,
\end{multline*}
where $\hat{}$ denotes the Fourier transform with respect to the
first $n$ variables. Since $\,b_{k,j}\in C_0^{n/2+\eps}(\R^{2n})\,$,
it follows that the functions $\,(1+|\eta|)^{n/2+\eps}\,\hat
b_{k,j}(\eta,\xi)\,$ belong to $\,L_2(\R^{2n})\,$ and their
$L_2$-norms are estimated by constants depending on $D$ and the
H\"older norms of $\si$ and $\psi_z\,$. Substituting this
representation in the corresponding oscillatory integral, we see
that $\,R_k=R^*_{k,1}R_{k,2}\,$ where the $R_{k,j}$ are
Hilbert-Schmidt operators acting from $\,L_2(M)\,$ to
$\,L_2(\R^{3n}\times[0,1]^2)\,$.

Since $\,\QC_1(\psi_0(\si))\in\Sf_2\,$ and $\,R_{k,j}\in\Sf_2\,$, we
have $\,\QC_{1/2}(\varf_0(\si))\in\Sf_1\,$. The Hilbert-Schmidt
norms of $\,R_{k,j}\,$ are bounded by $\,C_{k,j}\,h^{1-n}\,$ where
$\,C_{k,j}\,$ are constants depending only on $\delta$, $M$, $D$,
and the H\"older norms of $\si$ and $\varf$. Therefore
Theorem~\ref{T7} implies that
$$
\Tr\varf_0(\QC_{1/2}(\si))\ \leq\ (2\pi h)^{-n}\int_{T^*M}
\varf_0(\si)\,\dr x\,\dr\xi\;+\;C\,h^{1-n}\,,
$$
where $C$ is a constant depending on the same parameters. If
$\,\QC_{1/2}(\si)\in\Sf_1\,$ and $\,\varf(0)=0\,$ then the same
estimate holds for the function $\,\varf\,$.
\end{example}


\begin{thebibliography}{MMM}

\bibitem[1]{B1} F. Berezin.
{\it Convex functions of operators.} Mat. Sb. {\bf 88} (1972),
268--276 (Russian).

\bibitem[2]{B2} F. Berezin.
{\it Covariant and contravariant symbols of operators.} Izv. Akad.
Nauk SSSR Ser. Mat. {\bf 36} (1972), 1134--1167 (Russian).

\bibitem[3]{DS} M. Dimassi and J. Sj\"ostrand. {\it Spectral
Asymptotics in the Semiclassical Limit.} Cambridge University
Press, 1999.

\bibitem[4]{ELSS}
W.D. Evans, R.T. Lewis, H. Siedentop and J.P. Solovej {\it
Counting eigenvalues using coherent states with an application to
Dirac and Schr\"odinger operators in the semiclassical limit},
Ark. Mat. {\bf 34} (1996), 265--283.

\bibitem[5]{ET}
I. Ekeland and R. Temam.  {\it Convex analysis and variational
problems.} Studies in Mathematics and Its Applications, vol. 1.
North-Holland Publishing Company and American Elseveir Publishing
Company, 1976. Russian translation: ``Mir'', 1979.

\bibitem[6]{G}
J.R. Giles. {\it Convex analysis with application in the
differentiation of convex functions.} Research Notes in
Mathematics, vol. 58. Pitman, 1982.

\bibitem[7]{H}
L. H\"ormander. {\it The Analysis of Linear Partial Differential
Operators III.} Springer--Verlag, 1985. Russian translation:
``Mir'', 1987.

\bibitem[8]{LS1}
A. Laptev and Yu. Safarov. {\it A generalization of the
Berezin--Lieb inequality for convex functions,} Amer. Math. Soc.
Transl. Ser. 2 {\bf 175} (1996), 69--80.

\bibitem[9]{LS2}
A. Laptev and Yu. Safarov. {\it Szeg\"o type limit theorems,} J.
Funct. Anal. {\bf 138} (1996), 544--559.

\bibitem[10]{Sa}
Y. Safarov. {\it Pseudodifferential operators and linear
connections.} Proc. London Math. Soc., {\bf 74} (1997), 379--417.

\bibitem[11]{SV}
Yu. Safarov and D. Vassiliev. {\it The asymptotic distribution of
eigenvalues of partial differential operators.} American
Mathematical Society, 1996.

\bibitem[12]{T}
D. Tataru. {\it On the Fefferman--Phong inequality and related
problems.} Comm. Partial Differential Equations {\bf 27} (2002),
no. 11-12, 2101--2138.

\bibitem[13]{LM}
N. Lerner and Y. Morimoto. {\it On the Fefferman-Phong inequality
and a Wiener-type algebra of pseudodifferential operators.} Preprint
2005,
\newline
http://www.perso.univ-rennes1.fr/nicolas.lerner/FPI.2005.pdf
\newline
{\sl (The reference was added after the paper had been sent to
print, and therefore does not appear in the published text.)}


\end{thebibliography}
\end{document}